# Quotient groups of IA-automorphisms of free metabelian groups

C.E. Kofinas and A.I. Papistas

C.E. Kofinas, Department of Mathematics, University of the Aegean, Karlovassi, 832 00 Samos, Greece. *e-mail:* kkofinas@aegean.gr

A.I. Papistas, Department of Mathematics, Aristotle University of Thessaloniki, 541 24 Thessaloniki, Greece. *e-mail:* apapist@math.auth.gr

**Abstract**

For a positive integer $n$, with $n \geq 2$, let $M_n$ be a free metabelian group of rank $n$. For $c \in \mathbb{N}$, let $\gamma_c(M_n)$ be the $c$-th term of the lower central series of $M_n$. For $c \geq 2$, let $\mathrm{I}_c\mathrm{A}(M_n)$ be the subgroup of $\mathrm{Aut}(M_n)$ consisting of all automorphisms inducing the identity mapping on $M_n/\gamma_c(M_n)$. In this paper, we study the quotient groups $\mathcal{L}^c(\mathrm{IA}(M_n)) = \mathrm{I}_c\mathrm{A}(M_n)/\mathrm{I}_{c+1}\mathrm{A}(M_n)$ for all $n$ and $c$. For $c \geq 2$, we show $\mathrm{I}_c\mathrm{A}(M_2)) = \gamma_{c-1}(\mathrm{IA}(M_2))$. For $n = 3$, we show $\gamma_3(\mathrm{IA}(M_3)) \neq \mathrm{I}_4\mathrm{A}(M_3)$ and so, the Andreadakis' conjecture (for a free metabelian group) is not valid for $n = 3$ and $c = 3$. For $n \geq 4$ and $c \geq 3$, we prove that $\mathcal{L}^{c+1}(\mathrm{IA}(M_n)) = \gamma_c(\mathrm{IA}(M_n))\mathrm{I}_{c+1}\mathrm{A}(M_n)/\mathrm{I}_{c+1}\mathrm{A}(M_n)$.

## 1 Introduction

Let $G$ be a group. For subgroups $X$ and $Y$ of $G$, we write $[X,Y]$ for the subgroup of $G$ generated by all commutators $[x,y] = x^{-1}y^{-1}xy$ with $x \in X$ and $y \in Y$. For a positive integer $c$, let $\gamma_c(G)$ be the $c$-th term of the lower central series of $G$. We point out that $\gamma_2(G) = [G,G] = G'$; that is, the derived group of $G$, and $G'' = (G')'$. For a positive integer $c$, with $c \geq 3$ and $g_1, \ldots, g_c \in G$, we write $[g_1, \ldots, g_c] = [[g_1, \ldots, g_{c-1}], g_c]$. In particular, for $x, y \in G$ and a non-negative integer $m$, we write $[x, {}_m y] = [x, y, \ldots, y]$ with $m$ factors of $y$, and $[x, {}_0 y] = x$. We denote by $\mathrm{Aut}(G)$ the automorphism group of $G$. For $c \geq 2$, the natural epimorphism from $G$ onto $G/\gamma_c(G)$ induces a group homomorphism, say $\pi_{c,G}$, from $\mathrm{Aut}(G)$ into $\mathrm{Aut}(G/\gamma_c(G))$. Let $\mathrm{I}_c\mathrm{A}(G) = \mathrm{Ker}\pi_{c,G}$. That is, $\mathrm{I}_c\mathrm{A}(G)$ is the normal subgroup of $\mathrm{Aut}(G)$ consisting of all automorphisms of $G$ which induce the identity mapping on $G/\gamma_c(G)$. For $c = 2$, we write $\mathrm{I}_2\mathrm{A}(G) = \mathrm{IA}(G)$. The elements of $\mathrm{IA}(G)$ are called IA-automorphisms of $G$. It is shown in [1] that, for all $t, s \geq 2$, $[\mathrm{I}_t\mathrm{A}(G), \mathrm{I}_s\mathrm{A}(G)] \subseteq \mathrm{I}_{t+s-1}\mathrm{A}(G)$. By an inductive argument on $c$, we get $\gamma_c(\mathrm{IA}(G)) \subseteq \mathrm{I}_{c+1}\mathrm{A}(G)$ for all $c \geq 1$. For $c \geq 2$, we write $\mathcal{L}^c(\mathrm{IA}(G)) = \mathrm{I}_c\mathrm{A}(G)/\mathrm{I}_{c+1}\mathrm{A}(G)$ and for any subgroup $H$ of $\mathrm{IA}(G)$, let $\mathrm{gr}_{c-1}(H) = \gamma_{c-1}(H)/\gamma_c(H)$. In general, it is a very difficult problem to determine whether or not $\gamma_c(\mathrm{IA}(G)) = \mathrm{I}_{c+1}\mathrm{A}(G)$ for all $c \geq 2$. Moreover, it is quite difficult to study the structure of the quotient groups $\mathrm{gr}_c(\mathrm{IA}(G))$ for all $c \geq 1$.

For a positive integer $n$, with $n \geq 2$, let $F_n$ be a free group of rank $n$ with a free generating set $\{f_1, \ldots, f_n\}$. For $G = F_n$, the above problem is known as the Andreadakis' conjecture. It has been proved in [1, Theorem 6.1] that $\gamma_c(\mathrm{IA}(F_2)) = \mathrm{I}_{c+1}\mathrm{A}(F_2)$ for all $c \geq 2$, and that the first three terms of this series coincide in the case $n = 3$ (see, [1, Theorem 6.2]). For $n \geq 2$, $\gamma_2(\mathrm{IA}(F_n)) = \mathrm{I}_3\mathrm{A}(F_n)$ (see [15, Theorem 1.1]) and, for $n \geq 3$, $\gamma_3(\mathrm{IA}(F_n)) = \mathrm{I}_4\mathrm{A}(F_n)$ (see [17, Theorem 1], [15, Corollary 1.4]). For $n \geq 2$, let $M_n = F_n/F_n''$, that is, $M_n$ is a free metabelian group of rank $n$. The set $\{x_1, \ldots, x_n\}$, where $x_i = f_iF_n''$ and $i \in \{1, \ldots, n\}$, is a free generating set of $M_n$. The natural group epimorphism from $F_n$ onto $M_n$ induces a group homomorphism $\alpha_n$ from $\mathrm{Aut}(F_n)$ into $\mathrm{Aut}(M_n)$ and let $T_n = \mathrm{Im}\alpha_n$. Elements of $T_n$ are called tame automorphisms of $M_n$; otherwise the elements are called non-tame. In a series of papers ([1], [2], [7], [4], [5], [16]), it has been shown that $\alpha_n$, with $n \neq 3$, is onto and $\alpha_3$ is not onto.

In this paper, we investigate the quotient groups $\mathcal{L}^c(\mathrm{IA}(M_n))$, with $c \geq 2$. Our main purpose in this paper is to show the following result.

**Theorem 1** *Let $M_n$ be a free metabelian group of finite rank $n$, with $n \geq 2$. Then,*

1. *For $c \geq 1$, $\gamma_c(\mathrm{IA}(M_2)) = \mathrm{I}_{c+1}\mathrm{A}(M_2))$. In particular, for $c \geq 2$, $\mathcal{L}^c(\mathrm{IA}(M_2)) = \gamma_{c-1}(\mathrm{IA}(M_2))/\gamma_c(\mathrm{IA}(M_2))$.*

2. *$\mathcal{L}^4(\mathrm{IA}(M_3)) \neq \gamma_3(\mathrm{IA}(M_3))\mathrm{I}_5\mathrm{A}(M_3)/\mathrm{I}_5\mathrm{A}(M_3)$ and so, $\gamma_3(\mathrm{IA}(M_3)) \neq \mathrm{I}_4\mathrm{A}(M_3)$.*

3. *For $n \geq 4$ and $c \geq 3$, $\mathcal{L}^{c+1}(\mathrm{IA}(M_n)) = \gamma_c(\mathrm{IA}(M_n))\mathrm{I}_{c+1}\mathrm{A}(M_n)/\mathrm{I}_{c+1}\mathrm{A}(M_n)$.*

Let $\mu$ be the IA-endomorphism of $M_3$ satisfying the conditions $\mu(x_1) = x_1[x_1^{-1}, [x_1, [x_2, x_3]]]$ and $\mu(x_j) = x_j$, $j = 2, 3$. It is shown in [7, Proof of Theorem 8 and Theorem 2] that the IA-endomorphism $\mu$ of $M_3$ is a non-tame automorphism of $M_3$. Thus, $T_3 \neq \mathrm{Aut}(M_3)$. By a result of Bachmuth and Mochizuki [2, Theorem], $\mathrm{Aut}(M_3)$ is not finitely generated. By Theorem 1 (1), the Andreadakis' conjecture is valid for a free metabelian group of rank 2. Theorem 1 (2) gives information concerning the Andreadakis' conjecture (for a free metabelian group of rank 3); it is not valid for $n = 3$ and $c = 3$. For its proof, we first prove (Theorem 2) that, for $n \geq 2$ and $r \in \mathbb{N}$, $\gamma_r(H)(\mathrm{I}_{r+2}\mathrm{A}(G_n))/\mathrm{I}_{r+2}\mathrm{A}(G_n) = \gamma_r(\mathrm{IA}(G_n))(\mathrm{I}_{r+2}\mathrm{A}(G_n))/\mathrm{I}_{r+2}\mathrm{A}(G_n)$ for a certain subgroup $H$ of the automorphism group of a free non-abelian polynilpotent group $G_n$ of rank $n$. Next, by applying a result of Chein [7, Proof of Theorem 8], we obtain the required result. As we point out in Remark 1, the subgroup of $\mathrm{Aut}(M_3)$ generated by $T_3$ and $\mu$ is dense in $\mathrm{Aut}(M_3)$ with respect to the formal power series topology on $\mathrm{End}(M_3)$.

For $n \geq 4$ and $c \geq 2$, let $\rho_{n,c}$ be the natural epimorphism from $\mathrm{gr}_c(\mathrm{IA}(M_n))$ onto $\gamma_c(\mathrm{IA}(M_n))\mathrm{I}_{c+2}\mathrm{A}(M_n)/\mathrm{I}_{c+2}\mathrm{A}(M_n)$. That is, for all $\phi \in \gamma_c(\mathrm{IA}(M_n))$, $\rho_{n,c}(\phi\gamma_{c+1}(\mathrm{IA}(M_n))) = \phi\mathrm{I}_{c+2}\mathrm{A}(M_n)$. Since $\alpha_n$ is onto and $\gamma_c(\mathrm{IA}(F_n)) = \mathrm{I}_{c+1}\mathrm{A}(F_n)$ with $c = 2, 3$, we have $\gamma_c(\mathrm{IA}(M_n)) = \mathrm{I}_{c+1}\mathrm{A}(M_n)$ for $c = 2, 3$. By Theorem 1 (3), for all $c \geq 3$, $\rho_{n,c}$ is an epimorphism from $\mathrm{gr}_c(\mathrm{IA}(M_n))$ onto $\mathcal{L}^{c+1}(\mathrm{IA}(M_n))$. We point out that if $\rho_{n,c}$ is one-to-one for all $c \geq 3$, then $\gamma_c(\mathrm{IA}(M_n)) = \mathrm{I}_{c+1}\mathrm{A}(M_n)$ for all $c \geq 3$.

In the next few lines, we give an outline of the proof of Theorem 1 (3): In Section 4.3.1, for $c \geq 2$, we give on $\mathcal{L}^c(\mathrm{IA}(M_n))$ a structure of a $\mathbb{Z}\mathrm{GL}_n(\mathbb{Z})$-module in a natural way. For $c \in \mathbb{N}$, we consider the $\mathbb{Q}$-vector space $\mathrm{gr}_{c,\mathbb{Q}}(M_n)^{\oplus n} = \mathrm{gr}_{c,\mathbb{Q}}(M_n) \oplus \cdots \oplus \mathrm{gr}_{c,\mathbb{Q}}(M_n)$ with $n$ copies of $\mathrm{gr}_{c,\mathbb{Q}}(M_n)$. The general linear group $\mathrm{GL}_n(\mathbb{Q})$ acts on $\mathrm{gr}_{c,\mathbb{Q}}(M_n)^{\oplus n}$ as follows: For $g \in \mathrm{GL}_n(\mathbb{Q})$ and $(u_1, \ldots, u_n) \in \mathrm{gr}_{c,\mathbb{Q}}(M_n)^{\oplus n}$, $g * (u_1, \ldots, gu_n) = (gu_1, \ldots, gu_n)g^{-1}$, where $gu_i$ means the canonical action of $g$ on $\mathrm{gr}_{c,\mathbb{Q}}(M_n)$ and the multiplication of a $1 \times n$ and an $n \times n$ matrix. It follows from a result of Bryant and Drensky [6, Section 2, Proposition 3.5 (i)] that $\mathrm{gr}_{c,\mathbb{Q}}(M_n)^{\oplus n}$ is a (rational) $\mathbb{Q}\mathrm{GL}_n(\mathbb{Q})$-module. In fact, the module $\mathrm{gr}_{c,\mathbb{Q}}(M_n)^{\oplus n}$ is decomposed into three irreducible $\mathbb{Q}\mathrm{GL}_n(\mathbb{Q})$-modules, $P_c$, $Q_c$ and $R_c$. By means of a monomorphism $\overline{\chi}_c$, we embed $\mathcal{L}^c(\mathrm{IA}(M_n))$ into $\mathrm{gr}_{c,\mathbb{Q}}(M_n)^{\oplus n}$ as a $\mathbb{Z}\mathrm{GL}_n(\mathbb{Z})$-module. For $c \geq 3$, we observe that certain generators of $P_c$ and $Q_c$ may be regarded elements of $\mathcal{L}^c(\mathrm{IA}(M_n))$. By using a result of Bachmuth [3, Theorem 3 and Proof of Lemma 7], we prove that these certain IA-automorphisms of $M_n$ generate $\mathcal{L}^c(\mathrm{IA}(M_n))$ as a $\mathbb{Z}\mathrm{GL}_n(\mathbb{Z})$-module (Proposition 2). Finally, we show that these generators of $\mathcal{L}^c(\mathrm{IA}(M_n))$ belong to $\gamma_{c-1}(\mathrm{IA}(M_n))$ (Lemma 5) and so, we obtain the desired result. As a consequence, we have $\mathbb{Q} \otimes \overline{\chi}_c(\mathcal{L}^c(\mathrm{IA}(M_n))) = P_c \oplus Q_c$ as vector spaces over $\mathbb{Q}$ (Corollary 2).

# 2 Preliminaries

A group $G$ is poly-nilpotent if there exists in $G$ a series of subgroups: $G = N_0 > N_1 > \cdots > N_m = \{1\}$, where each $N_\alpha$ is normal in $N_{\alpha-1}$ and $N_{\alpha-1}/N_\alpha$ is nilpotent. If $N_{\alpha-1}/N_\alpha$ has class $c_\alpha - 1$, then $\gamma_{c_\alpha}(N_{\alpha-1}) \leq N_\alpha$ and we note that $c_\alpha > 1$. Writing $\gamma_k\gamma_h(H) = \gamma_k(\gamma_h(H))$, where $H$ is any subgroup of $G$, we have $\gamma_{c_m} \cdots \gamma_{c_1}(G) = \{1\}$. For a positive integer $n \geq 2$, write $N = \gamma_{c_m} \cdots \gamma_{c_1}(F_n)$, with $c_1 \geq 3$, and let $G_n = F_n/N$. Thus, $G_n$ is a free non-abelian poly-nilpotent group of finite rank $n$ (of class row $(c_1 - 1, \ldots, c_m - 1)$). For $i \in [n] = \{1, \ldots, n\}$, we write $y_i = f_iN$ and so, the set $\{y_1, \ldots, y_n\}$ is a free generating set of $G_n$. Note that $G_n$ is residually nilpotent (see [9]) and each quotient group $\mathrm{gr}_c(G_n)$ is a free abelian group of finite rank (see [18]). For all $c \geq 2$, let $G_{n,c} = G_n/\gamma_c(G_n)$. For $i \in [n]$, let $y_{i,c} = y_i\gamma_c(G_n)$. So, $\{y_{1,c}, \ldots, y_{n,c}\}$ is a free generating set of $G_{n,c}$. Since $G_{n,c+1}/\mathrm{gr}_c(G_n) \cong G_{n,c}$ in a natural way and $\mathrm{gr}_c(G_n)$ is a fully invariant subgroup of $G_{n,c+1}$, the natural group epimorphism from $G_{n,c+1}$ onto $G_{n,c}$ induces a group epimorphism $\psi_{c+1,G}$ from $\mathrm{Aut}(G_{n,c+1})$ onto $\mathrm{Aut}(G_{n,c})$. Recall that $\pi_{c,G}$ is the induced homomorphism from $\mathrm{Aut}(G_n)$ into $\mathrm{Aut}(G_{n,c})$. We define $A_{c+1}(G_n) = \mathrm{Ker}\psi_{c+1,G}$ and $A^*_{c+1}(G_n) = \mathrm{Im}\pi_{c+1,G} \cap \mathrm{Ker}\psi_{c+1,G}$. For $t \in \{2, \ldots, c\}$, the natural group epimorphism from $G_{n,c}$ onto $G_{n,c}/\gamma_t(G_{n,c})$ induces a group homomorphism $\sigma_{c,t} : \mathrm{Aut}(G_{n,c}) \to \mathrm{Aut}(G_{n,c}/\gamma_t(G_{n,c}))$. We write $\mathrm{I}_t\mathrm{A}(G_{n,c}) = \mathrm{Ker}\sigma_{c,t}$ and, for $t = 2$, $\mathrm{IA}(G_{n,c}) = \mathrm{I}_2\mathrm{A}(G_{n,c})$. Thus, for $c \geq 2$, $A_{c+1}(G_n) = \mathrm{I}_c\mathrm{A}(G_{n,c+1})$ and $A^*_{c+1}(G_n) = \mathrm{Im}\pi_{c+1,G} \cap \mathrm{I}_c\mathrm{A}(G_{n,c+1})$. By definitions, for all positive integers $n, c$, with $n, c \geq 2$, $A_{c+1}(G_n)$ is isomorphic as an abelian group to $\mathrm{gr}_c(G_n)^{\oplus n}$, a direct sum of $n$ copies of $\mathrm{gr}_c(G_n)$, and $A^*_{c+1}(G_n) \cong \mathrm{I}_c\mathrm{A}(G_n)/\mathrm{I}_{c+1}\mathrm{A}(G_n)$ as abelian groups in a natural way. Since, for all $c \geq 1$, $\mathrm{gr}_c(G_n))$ is a free abelian group, we have $A_{c+1}(G_n)$ is a free abelian group of finite rank for all $n$ and $c$, with $n \geq 2$ and $c \geq 2$, and so, $A^*_{c+1}(G_n)$ is a free abelian group of finite rank for all $n$ and $c$, with $n \geq 2$ and $c \geq 2$. It is important to point out that all the above definitions and constructions concerning $G_n$ are valid for any free group $F_n$ of rank $n \geq 2$ as well.

For $n, c \geq 2$, the set $\mathcal{M}_{n,c} = \{[x_{i_1}, x_{i_2}, \ldots, x_{i_c}]\gamma_{c+1}(M_n) : i_1 > i_2 \leq i_3 \leq \cdots \leq i_c; i_1, \ldots, i_c \in [n]\}$ is a $\mathbb{Z}$-basis of $\mathrm{gr}_c(M_n)$ and the rank of $\mathrm{gr}_c(M_n)$ is $(c-1)\binom{n+c-2}{n-2}$ (see, for example, [3]) and so, the rank of $A_{c+1}(M_n)$ is $n(c-1)\binom{n+c-2}{n-2}$. It follows from a result of Andreadakis [1, Theorem 5.1] (see, also, [3, Theorem 1]) that the rank of $A^*_3(M_n)$ is $n\binom{n}{2}$. By a result of Bachmuth [3, Theorem 2], we have the rank of $A^*_{c+1}(M_n)$, with $c \geq 3$, is $n(c-1)\binom{n+c-2}{n-2} - \binom{n+c-2}{n-1}$.

For $c \geq 2$, let $\mathcal{L}^c(\mathrm{IA}(G_n)) = \mathrm{I}_c\mathrm{A}(G_n)/\mathrm{I}_{c+1}\mathrm{A}(G_n)$ and we write the product in each of $\mathcal{L}^c(\mathrm{IA}(G_n))$ additively. Namely, for any $\alpha, \beta \in \mathrm{I}_c\mathrm{A}(G_n)$, if we denote their coset classes modulo $\mathrm{I}_{c+1}\mathrm{A}(G_n)$ by $\overline{\alpha}$ and $\overline{\beta}$, respectively, then $\overline{\alpha} + \overline{\beta} = \overline{\alpha\beta}$. Form the (restricted) direct sum of the abelian groups $\mathcal{L}^c(\mathrm{IA}(G_n))$, denoted by $\mathcal{L}(\mathrm{IA}(G_n)) = \bigoplus_{c \geq 2} \mathcal{L}^c(\mathrm{IA}(G_n))$. It has the structure of a graded Lie ring with $\mathcal{L}^c(\mathrm{IA}(G_n))$ as component of degree $c-1$ in the grading and Lie multiplication given by $[\phi \mathrm{I}_{j+1}\mathrm{A}(G_n), \psi \mathrm{I}_{\kappa+1}\mathrm{A}(G_n)] = [\phi, \psi]\mathrm{I}_{j+\kappa}\mathrm{A}(G_n)$ $(j, \kappa \geq 2)$. For any subgroup $H$ of $\mathrm{Aut}(G_n)$ and a positive integer $r$, with $r \geq 2$, let $H_r = H \cap \mathrm{I}_r\mathrm{A}(G_n)$. That is, $H_r$ consists of all elements of $H$ which induce the identity map on $G_{n,r}$. For $q \in \mathbb{N}$, let $\mathcal{L}^q_1(H) = \gamma_q(H)(\mathrm{I}_{q+2}\mathrm{A}(G_n))/\mathrm{I}_{q+2}\mathrm{A}(G_n)$ and, for $q \geq 2$, $\mathcal{L}^q_2(H) = H_q(\mathrm{I}_{q+1}\mathrm{A}(G_n))/\mathrm{I}_{q+1}\mathrm{A}(G_n)$. Form the (restricted) direct sums of abelian groups $\mathcal{L}^q_1(H)$ and $\mathcal{L}^q_2(H)$, denoted by $\mathcal{L}_1(H)$ and $\mathcal{L}_2(H)$, respectively. We point out that, for all $r \in \mathbb{N}$, $\mathcal{L}^r_1(H) \subseteq \mathcal{L}^{r+1}_2(H) \subseteq \mathcal{L}^{r+1}(\mathrm{IA}(G_n))$. Clearly, both $\mathcal{L}_1(H)$ and $\mathcal{L}_2(H)$ are Lie subrings of $\mathcal{L}(\mathrm{IA}(G_n))$ and $\mathcal{L}_1(H) \subseteq \mathcal{L}_2(H)$. Notice that if $H = \mathrm{IA}(G_n)$, then $\mathcal{L}_2(\mathrm{IA}(G_n)) = \mathcal{L}(\mathrm{IA}(G_n))$. In [8], the formal power series topology for endomorphisms of relatively free groups is introduced, simulating the approach of [6]. Let $H_1, H_2 \leq \mathrm{Aut}(G_n)$ subject to $T_{n,G} \subseteq H_1 \subseteq H_2$, where $T_{n,G}$ is the subgroup of $\mathrm{Aut}(G_n)$ consisting of those automorphisms of $G_n$ which are induced by automorphisms of $F_n$. As observed in [8, Proposition 2.1 (ii)], $H_1$ is dense in $H_2$ with respect to the formal power series topology on $\mathrm{End}(G_n)$ if and only if $\mathcal{L}_2(H_1) = \mathcal{L}_2(H_2)$.

For any group $G$ and a positive integer $c$, we write $\mathrm{gr}(G) = \bigoplus_{c \geq 1} \mathrm{gr}_c(G)$ for the (restricted) direct sum of the abelian groups $\mathrm{gr}_c(G)$. It is well known that $\mathrm{gr}(G)$ has the structure of a Lie ring by defining a Lie multiplication $[a\gamma_{r+1}(G), b\gamma_{s+1}(G)] = [a, b]\gamma_{r+s+1}(G)$, where $a\gamma_{r+1}(G)$ and $b\gamma_{s+1}(G)$ are the images of the elements $a \in \gamma_r(G)$ and $b \in \gamma_s(G)$ in the quotient groups $\mathrm{gr}_r(G)$ and $\mathrm{gr}_s(G)$, respectively, and $[a, b]\gamma_{r+s+1}(G)$ is the image of the group commutator $[a, b]$ in the quotient group $\mathrm{gr}_{r+s}(G)$. Multiplication is then extended to $\mathrm{gr}(G)$ by linearity.

# 3 An observation on free poly-nilpotent groups

Throughout this section, for a positive integer $n \geq 2$, we write $G_n = F_n/N$, with $N = \gamma_{c_m} \cdots \gamma_{c_1}(F_n)$ and $c_1 \geq 3$, for a free non-abelian poly-nilpotent group of rank $n$. The natural epimorphism from $F_n$ onto $G_n$ induces a group homomorphism $\rho_n$ from $\mathrm{Aut}(F_n)$ into $\mathrm{Aut}(G_n)$. For $\beta \in \mathrm{Aut}(F_n)$, $\rho_n(\beta)(gN) = \beta(g)N$ for all $g \in F_n$. Write $T_{n,G} = \mathrm{Im}\rho_n$, that is, $T_{n,G}$ consists of all automorphisms of $G_n$ induced by automorphisms of $F_n$. The elements of $T_{n,G}$ are called tame automorphisms of $G_n$. For $n \geq 2$, we denote $H_{n,G} = T_{n,G} \cap \mathrm{IA}(G_n)$. For a positive integer $r$, with $r \geq 2$, let $H_{n,G,r} = T_{n,G} \cap \mathrm{I}_r\mathrm{A}(G_n)$. That is, $H_{n,G,r}$ consists of all tame automorphisms of $G_n$ which induce the identity mapping on $G_{n,r}$. Note that $H_{n,G,r} = H_{n,G} \cap \mathrm{I}_r\mathrm{A}(G_n)$ for $r \geq 2$. In the case $G_n = M_n$, we simply write $T_{n,G} = T_n$, $H_{n,G} = H_n$ and $H_{n,G,r} = H_{n,r}$. Since $G_n/G_n' \cong F_n/F_n'$, $\rho_n$ induces a group homomorphism $\rho_n' : \mathrm{IA}(F_n) \to \mathrm{IA}(G_n)$.

**Lemma 1** *With the previous notations, for all $n$, with $n \geq 2$, $H_{n,G} = \rho_n'(\mathrm{IA}(F_n))$ and $\gamma_2(H_{n,G}) = H_{n,G,3} = \rho_n'(\mathrm{I_3A}(F_n))$.*

*Proof.* Let $\alpha \in H_{n,G}$. Then, there exists $\beta \in \mathrm{Aut}(F_n)$ such that $\rho_n(\beta) = \alpha$. Since $\alpha \in \mathrm{IA}(G_n)$ and $G_n/G_n' \cong F_n/F_n'$, we obtain $\beta \in \mathrm{IA}(F_n)$. Thus, $\alpha = \rho_n'(\beta)$ and so, $H_{n,G} \subseteq \rho_n'(\mathrm{IA}(F_n))$. On the other hand, it is clear that $\rho_n'(\beta) \in H_{n,G}$ for all $\beta \in \mathrm{IA}(F_n)$. Therefore, $H_{n,G} = \rho_n'(\mathrm{IA}(F_n))$ and, for all $r \geq 2$, $H_{n,G,r} = \rho_n'(\mathrm{IA}(F_n)) \cap \mathrm{I}_r\mathrm{A}(G_n)$. Since $\mathrm{I_3A}(F_n) = \gamma_2(\mathrm{IA}(F_n))$ and $H_{n,G} = \rho_n'(\mathrm{IA}(F_n))$, we get

$$\gamma_2(H_{n,G}) = [H_{n,G}, H_{n,G}] = [\rho_n'(\mathrm{IA}(F_n)), \rho_n'(\mathrm{IA}(F_n))] = \rho_n'(\gamma_2(\mathrm{IA}(F_n))) = \rho_n'(\mathrm{I_3A}(F_n)).$$

Next, we claim that $\gamma_2(H_{n,G}) = H_{n,G,3}$. Since $\gamma_2(H_{n,G}) = \rho_n'(\mathrm{I_3A}(F_n))$, it is enough to show that $H_{n,G,3} = \rho_n'(\mathrm{I_3A}(F_n))$. Indeed, since $\gamma_3(G_n) = \gamma_3(F_n)/N$, we have $\rho_n'(\mathrm{I_3A}(F_n)) \subseteq H_{n,G,3}$. On the other hand, let $\alpha \in H_{n,G,3}$. Then, there exists $\beta \in \mathrm{IA}(F_n)$ such that $\rho_n'(\beta) = \alpha$. That is, $\rho_n'(\beta)(y_i) = \beta(f_i)N = \alpha(y_i)$ for all $i \in [n]$. Since $\alpha \in \mathrm{I_3A}(G_n)$, we get $\alpha(y_i) = y_i v_i = f_i u_i N$, with $u_i \in \gamma_3(F_n)$ for all $i \in [n]$. Hence, $\beta \in \mathrm{I_3A}(F_n)$ and so, $\alpha \in \rho_n'(\mathrm{I_3A}(F_n))$. Therefore, $H_{n,G,3} \subseteq \rho_n'(\mathrm{I_3A}(F_n))$ and so, $H_{n,G,3} = \rho_n'(\mathrm{I_3A}(F_n))$. □

In the following result, we obtain information about $\mathcal{L}_1^1(H_{n,G})$ for all $n$, with $n \geq 2$.

**Lemma 2** *For any positive integer $n$, with $n \geq 2$, $\mathcal{L}_1^1(H_{n,G}) = \mathcal{L}^2(\mathrm{IA}(G_n))$. Moreover, $\mathcal{L}_1^1(H_{n,G})$ is a free abelian group of rank $\frac{n^2(n-1)}{2}$.*

*Proof.* It has been shown in ([1, Theorem 3.1], [3, Theorem 1]) that $A_3^*(F_n) = A_3(F_n)(\cong \mathrm{gr}_2(F_n)^{\oplus n})$ for all $n \geq 2$. Thus, $\mathrm{rank}(A_3(F_n)) = n\,\mathrm{rank}(\mathrm{gr}_2(F_n)) = \frac{n^2(n-1)}{2} = \mathrm{rank}(A_3^*(F_n)) = \mathrm{rank}(\mathcal{L}^2(\mathrm{IA}(F_n)))$. Since $G_n = F_n/N$ and $N \subseteq \gamma_3(F_n)$, we get $\mathrm{gr}_2(G_n) \cong \mathrm{gr}_2(F_n)$ as free abelian groups in a natural way and so, $A_3(G_n) \cong A_3(F_n)$ as free abelian groups. We point out that $G_{n,3} \cong F_{n,3}$ as groups in a natural way. Any automorphism of $F_{n,3}$ is tame, that is, $\pi_{3,F}(\mathrm{Aut}(F_n)) = \mathrm{Aut}(F_{n,3})$ (see, [1, Theorem 3.1]). Since $A_3(F_n) = A_3^*(F_n)$ and $G_{n,3} \cong F_{n,3}$ as groups, we have $A_3(G_n) = A_3^*(G_n)$. Therefore, $A_3(G_n) \cong \mathcal{L}^2(\mathrm{IA}(G_n))$ as free abelian groups and so, $\mathrm{rank}(\mathcal{L}^2(\mathrm{IA}(G_n))) = \frac{n^2(n-1)}{2}$.

Since, by Lemma 1, $\mathcal{L}_1^1(H_{n,G}) \subseteq \mathcal{L}^2(\mathrm{IA}(G_n))$, it is enough to show that $\mathcal{L}^2(\mathrm{IA}(G_n)) \subseteq \mathcal{L}_1^1(H_{n,G})$. We first calculate the rank of $\mathcal{L}_1^1(H_{n,G})$. Since, by Lemma 1, $\gamma_2(H_{n,G}) = H_{n,G,3}$, we have $\mathcal{L}_1^1(H_{n,G}) \cong H_{n,G}/\gamma_2(H_{n,G})$ as abelian groups. Therefore, $H_{n,G}/\gamma_2(H_{n,G})$ is a free abelian group. We claim that $H_{n,G}/\gamma_2(H_{n,G}) \cong \mathrm{gr}_2(G_n)^{\oplus n}$ as groups. We define a map $\omega : H_{n,G} \longrightarrow \mathrm{gr}_2(G_n)^{\oplus n}$ by $\omega(\alpha) = (y_1^{-1}\alpha(y_1)\gamma_3(G_n), \ldots, y_n^{-1}\alpha(y_n)\gamma_3(G_n))$ for all $\alpha \in H_{n,G}$. Note that $\alpha(y_i) = y_i u_i$, $u_i \in \gamma_2(G_n)$, $i = 1, \ldots, n$. Clearly, $\omega$ is a group homomorphism. Since $\gamma_2(H_{n,G}) = H_{n,G,3}$, we have $\mathrm{Ker}\omega = \gamma_2(H_{n,G})$. Thus, $\omega$ induces a group monomorphism $\overline{\omega} : H_{n,G}/\gamma_2(H_{n,G}) \to \mathrm{gr}_2(G_n)^{\oplus n}$. Next, we claim that $\overline{\omega}$ is onto. Let $(\overline{v}_1, \ldots, \overline{v}_n)$ be a non-trivial element of $\mathrm{gr}_2(G_n)^{\oplus n}$, where $\overline{v}_i = v_i\gamma_3(G_n)$, $v_i \in \gamma_2(G_n)$, $i \in [n]$. For $i \in [n]$, write $v_i = u_i N$, with $u_i \in \gamma_2(F_n)$. Let $\psi$ be the IA-automorphism of $F_{n,3}$ satisfying the conditions $\psi(f_i\gamma_3(F_n)) = f_i u_i \gamma_3(F_n)$, with $i \in [n]$. That is, $\psi \in A_3(F_n) = \mathrm{IA}(F_{n,3})$. Since $\pi_{3,F} : \mathrm{Aut}(F_n) \to \mathrm{Aut}(F_{n,3})$ is surjective and $F_{n,3}/F'_{n,3} \cong F_n/F'_n$, there exists $\beta \in \mathrm{IA}(F_n)$ such that $\pi_{3,F}(\beta) = \psi$. Let $\alpha = \rho'_n(\beta)$. Thus, $\alpha \in H_{n,G} = \rho'_n(\mathrm{IA}(F_n))$. Clearly, $\overline{\omega}(\alpha\gamma_2(H_{n,c})) = (\overline{v}_1, \ldots, \overline{v}_n)$ and so, $\overline{\omega}$ is onto. Therefore, $H_{n,G}/\gamma_2(H_{n,G}) \cong \mathrm{gr}_2(G_n)^{\oplus n}$. Since $\mathrm{gr}_2(G_n)^{\oplus n} \cong \mathcal{L}^2(\mathrm{IA}(G_n))$, we get $\mathrm{rank}(H_{n,G}/\gamma_2(H_{n,G})) = \mathrm{rank}(\mathcal{L}_1^1(H_{n,G})) = \mathrm{rank}(\mathcal{L}^2(\mathrm{IA}(G_n))) = \frac{1}{2}n^2(n-1)$. Since $\mathcal{L}_1^1(H_{n,G}) \subseteq \mathcal{L}^2(\mathrm{IA}(G_n))$ and $\mathcal{L}_1^1(\mathrm{IA}(G_n))$ is a free abelian group of the same rank, we obtain $\mathcal{L}_1^1(H_{n,G})$ has a finite index in $\mathcal{L}^2(\mathrm{IA}(G_n))$. Write $m = \frac{1}{2}n^2(n-1)$ and let $a_1, \ldots, a_m \in \mathrm{IA}(G_n)$ such that $\{a_1\mathrm{I}_3\mathrm{A}(G_n), \ldots, a_m\mathrm{I}_3\mathrm{A}(G_n)\}$ is a $\mathbb{Z}$-basis of $\mathcal{L}^2(\mathrm{IA}(G_n))$. Since $\pi_{3,F}$ : $\mathrm{Aut}(F_n) \to \mathrm{Aut}(F_{n,3})$ is surjective and $F_{n,3}/F'_{n,3} \cong F_n/F'_n$, we have $\pi_{3,F}$ induces an epimorphism $\pi'_{3,F} : \mathrm{IA}(F_n) \to \mathrm{IA}(F_{n,3})$. Since $A_3(G_n) \cong \mathrm{gr}_2(F_n)^{\oplus n}$, $A_3(G_n) = A_3^*(G_n) \cong \mathcal{L}^2(\mathrm{IA}(G_n))$, $\pi'_{3,F}$ is surjective and $G_{n,3}/G'_{n,3} \cong F_n/F'_n$, we may choose $a_1, \ldots, a_m$ to be elements of $H_{n,G}$. Therefore, $\mathcal{L}_1^1(H_{n,G}) = \mathcal{L}^2(\mathrm{IA}(G_n))$. □

Next, we show the main result in this section.

**Theorem 2** *For a positive integer n, with $n \geq 2$, let $G_n$ be a free non-abelian poly-nilpotent group of rank n. Then, for all integers n and r, with $n \geq 2$, $\mathcal{L}_1^r(H_{n,G}) = \mathcal{L}_1^r(\mathrm{IA}(G_n))$.*

*Proof.* It is enough to show that $\mathcal{L}_1(H_{n,G}) = \mathcal{L}_1(\mathrm{IA}(G_n))$. Moreover, since $\mathcal{L}_1(H_{n,G}) \subseteq \mathcal{L}_1(\mathrm{IA}(G_n))$, it is enough to show that $\mathcal{L}_1(\mathrm{IA}(G_n)) \subseteq \mathcal{L}_1(H_{n,G})$. Since $H_{n,G}/\gamma_2(H_{n,G})$ is a free abelian group of rank $m(= \frac{1}{2}n^2(n-1))$ (see proof of Lemma 2), we obtain $\mathrm{gr}(H_{n,G})$ is a finitely generated Lie ring (see, for example, [12, Lemma 3.2.3]). In fact, any generating set of $H_{n,G}/\gamma_2(H_{n,G})$ creates a generating set for $\mathrm{gr}(H_{n,G})$. Let $\mathcal{H}_G = \{b_1\gamma_2(H_{n,G}), \ldots, b_m\gamma_2(H_{n,G})\}$ be a free generating set of $H_{n,G}/\gamma_2(H_{n,G})$. For $r \in \mathbb{N}$, there is a natural group epimorphism from $\mathrm{gr}_r(H_{n,G})$ onto $\gamma_r(H_{n,G})/(\gamma_r(H_{n,G}) \cap \mathrm{I}_{r+2}\mathrm{A}(G_n))$. But $\gamma_r(H_{n,G})/(\gamma_r(H_{n,G}) \cap \mathrm{I}_{r+2}\mathrm{A}(G_n)) \cong \mathcal{L}_1^r(H_{n,G})$ as abelian groups in a natural way. Therefore, for all $r \geq 1$, there exists a natural group epimorphism $\zeta_r : \mathrm{gr}_r(H_{n,G}) \to \mathcal{L}_1^r(H_{n,G})$. Hence, the Lie ring $\mathcal{L}_1(H_{n,G})$ is generated by the $\zeta_1(b_1\gamma_2(H_{n,G})), \ldots, \zeta_1(b_m\gamma_2(H_{n,G}))$ (see, for example, [11, Lemma 9.4]). Since $\zeta_1(b_j\gamma_2(H_{n,G})) = b_j(\mathrm{I}_3\mathrm{A}(G_n))$ for $j \in \{1, \ldots, m\}$, we have $\mathcal{L}_1(H_{n,G})$ is generated by $b_1(\mathrm{I}_3\mathrm{A}(G_n)), \ldots, b_m(\mathrm{I}_3\mathrm{A}(G_n))$.

As in the proof of Lemma 2, we choose $b_1, \ldots, b_m \in H_{n,G}$ subject to the set $\{b_1(\mathrm{I}_3\mathrm{A}(G_n)), \ldots, b_m(\mathrm{I}_3\mathrm{A}(G_n))\}$ is a $\mathbb{Z}$-basis of $\mathcal{L}^2(\mathrm{IA}(G_n))$. For $r \in \mathbb{N}$, we claim that $\mathcal{L}_1^r(\mathrm{IA}(G_n))$ is generated as a group by the set $\mathcal{Z}_r = \{[b_{i_1}, \ldots, b_{i_r}]\mathrm{I}_{r+2}\mathrm{A}(G_n) : i_1, \ldots, i_r \in \{1, \ldots, m\}\}$. Clearly, we assume that $r \geq 2$. Let $h \in \gamma_r(\mathrm{IA}(G_n))$. Then, $h = \prod_{\mathrm{finite}}[g_1, \ldots, g_r]^{a_g}$, where $a_g \in \mathbb{Z}$, $g_1, \ldots, g_r \in \mathrm{IA}(G_n)$. Write each $g_j = \prod_{i=1}^m b_i^{n_{ij}} v_{ij}$, where $n_{ij} \in \mathbb{Z}$, $v_{ij} \in \mathrm{I}_3\mathrm{A}(G_n)$, $j = 1, \ldots, r$. Since, in any group $G$, $[ab, cd] = [a, d][a, d, b][b, d][a, c][a, c, bd][b, c][b, c, d]$ for all $a, b, c, d \in G$, $(\mathrm{I}_s\mathrm{A}(G_n))_{s\geq 2}$ is a central series of $\mathrm{IA}(G_n)$ and $\gamma_\mu(\mathrm{IA}(G_n)) \subseteq \mathrm{I}_{\mu+1}\mathrm{A}(G_n)$ for all $\mu$, we have

$$[g_1, \ldots, g_r]\mathrm{I}_{r+2}\mathrm{A}(G_n) = \prod_{\text{finite}} [b_{j_1}, \ldots, b_{j_r}]^{\pm 1}\mathrm{I}_{r+2}\mathrm{A}(G_n)$$

and so, $\mathcal{L}_1^r(\mathrm{IA}(G_n))$ is generated as a group by the set $\mathcal{Z}_r$. Thus, $\mathcal{L}_1(\mathrm{IA}(G_n))$ is generated as a Lie ring by the set $\{b_1(\mathrm{I}_3\mathrm{A}(G_n)), \ldots, b_m(\mathrm{I}_3\mathrm{A}(G_n))\}$. Since, by Lemma 2, $\mathcal{L}_1^1(H_{n,G}) = \mathcal{L}^2(\mathrm{IA}(G_n)) = \mathcal{L}_1^1(\mathrm{IA}(G_n))$, we have $\mathcal{L}_1(\mathrm{IA}(G_n)) \subseteq \mathcal{L}_1(H_{n,G})$ and so, we obtain the required result. □

# 4 The quotient groups $\mathcal{L}^c(\mathrm{IA}(M_n))$ for $n, c \geq 2$

## 4.1 The case $n = 2$

By combining results from [1, Theorem 6.1] and [2, Theorem 2], we have Theorem 1 for $n = 2$. In this section, we obtain the desired result in a more general way. In the next few lines, for any group $G$, we denote by $Z(G)$ the center of $G$. Let $\mathrm{Inn}(G)$ be the group of inner automorphisms of $G$. For a positive integer $c$, with $c \geq 2$, let $\mathrm{Inn}_c(G) = \mathrm{Inn}(G) \cap \mathrm{I}_c\mathrm{A}(G)$. That is, $\mathrm{Inn}_c(G)$ consists of all inner automorphisms of $G$ which induce the identity map on $G/\gamma_c(G)$. Our main claim in this section follows from the following result and a result of Bachmuth [2, Theorem 2].

**Lemma 3** *Let $G$ be a group. If $Z(G/\gamma_{c+1}(G)) = \mathrm{gr}_c(G)$ for all $c$, then, $\gamma_d(\mathrm{Inn}(G)) = \mathrm{Inn}_{d+1}(G)$ for all $d$.*

*Proof.* For $g \in G$, let $\tau_g \in \mathrm{Inn}_{d+1}(G)$ be the inner automorphism of $g$. That is, $\tau_g(x) = gxg^{-1}$ for all $x \in G$. Since $\tau_g \in \mathrm{Inn}_{d+1}(G)$, we have $[x, g^{-1}] \in \gamma_{d+1}(G)$ for all $x \in G$. Since $Z(G/\gamma_{d+1}(G)) = \mathrm{gr}_d(G)$, we get $g^{-1}\gamma_{c+1}(G) \in Z(G/\gamma_{d+1}(G)$ and so, $g \in \gamma_d(G)$. Write $g$ as a product of group commutators of the form $[g_1, \ldots, g_d]$ with $g_1, \ldots, g_d \in G$. Since $[\tau_x, \tau_y] = \tau_{[x,y]}$ for all $x, y \in G$, we get $\tau_g \in \gamma_d(\mathrm{Inn}(G))$ and so, $\mathrm{Inn}_{d+1}(G) \subseteq \gamma_d(\mathrm{Inn}(G))$. On the other hand, since $\gamma_d(\mathrm{Inn}(G)) \subseteq \gamma_d(\mathrm{IA}(G)) \subseteq \mathrm{I}_{d+1}\mathrm{A}(G)$, we obtain $\gamma_d(\mathrm{Inn}(G)) \subseteq \mathrm{Inn}_{d+1}(G)$ and so, $\gamma_d(\mathrm{Inn}(G)) = \mathrm{Inn}_{d+1}(G)$. □

Since $\mathrm{IA}(M_2) = \mathrm{Inn}(M_2)$ (see, [2, Theorem 2]) and $Z(M_2/\gamma_{c+1}(M_2)) = \mathrm{gr}_c(M_2)$ for all $c$ (see, for example, [14, Chapter 3, Section 6]), we have, by Lemma 3, the following result.

**Corollary 1** *For a positive integer $c$, $\gamma_c(\mathrm{IA}(M_2)) = \mathrm{I}_{c+1}\mathrm{A}(M_2)$.*

## 4.2 The case $n = 3$

Throughout this section, we write $H_3 = T_3 \cap \mathrm{IA}(M_3) = \rho_3'(\mathrm{IA}(F_3))$ (by Lemma 1) and, for all $r \in \mathbb{N}$, with $r \geq 2$, $H_{3,r} = T_3 \cap \mathrm{I}_r\mathrm{A}(M_3)$. Note that $H_3 = H_{3,2}$. It can be deduced from [7, Proof of Theorem 8] that $\mathcal{L}_2^4(H_3) \neq \mathcal{L}^4(\mathrm{IA}(M_3)) \cong A_5^*(M_3)$. By a result of Andreadakis [1, Theorem 6.2], we have

$$\gamma_3(H_3) = \gamma_3(\rho_3'(\mathrm{IA}(F_3))) = \rho_3'(\gamma_3(\mathrm{IA}(F_3))) = \rho_3'(\mathrm{I}_4\mathrm{A}(F_3)). \tag{3.1}$$

Clearly, $\gamma_3(H_3) \subseteq H_{3,4}$. Since $F_3/F_3' \cong M_3/M_3'$, it follows from the equation (3.1) that $H_{3,4} \subseteq \rho_3'(\mathrm{I}_4\mathrm{A}(F_3)) = \gamma_3(H_3)$. Therefore, $\gamma_3(H_3) = H_{3,4}$ and so,

$$\mathcal{L}_1^3(H_3) = \mathcal{L}_2^4(H_3). \tag{3.2}$$

By Theorem 2, we have $\mathcal{L}_1^3(H_3) = \mathcal{L}_1^3(\mathrm{IA}(M_3))$. By the equation (3.2) and since $\mathcal{L}_2^4(H_3) \neq \mathcal{L}^4(\mathrm{IA}(M_3))$, we obtain $\mathcal{L}_1^3(\mathrm{IA}(M_3)) \neq \mathcal{L}^4(\mathrm{IA}(M_3))$. Hence, $\gamma_3(\mathrm{IA}(M_3)) \neq \mathrm{I}_4\mathrm{A}(M_3)$.

From the aforementioned discussion, we deduce the following result.

**Proposition 1** *With the previous notations,* $\gamma_3(\mathrm{IA}(M_3)) \neq \mathrm{I}_4\mathrm{A}(M_3)$.

**Remark 1** Let $T_{3,\mu}$ be the subgroup of $\mathrm{Aut}(M_3)$ generated by $T_3$ and the IA-automorphism $\mu$ of $M_3$ satisfying the conditions $\mu(x_1) = x_1[x_1^{-1}, [x_1, [x_2, x_3]]]$ and $\mu(x_j) = x_j$, $j = 2, 3$. As pointed out in [8, Example 4.2], $T_3$ is not dense in $\mathrm{Aut}(M_3)$. By a result of Chein [7, Theorem 12], $\mathcal{L}_2(T_{3,\mu}) = \mathcal{L}_2(\mathrm{Aut}(M_3))$ and so, $T_{3,\mu}$ is dense in $\mathrm{Aut}(M_3)$.

## 4.3 The case $n \geq 4$

The main purpose in this section is to give a proof of Theorem 1 (3). For a group $G$, we write $\mathbb{Z}G$ for the integral group ring of $G$. For elements $x$ and $y$ of a group $G$, $y^x$ denotes the conjugate $x^{-1}yx$. Since $M_n'$ is abelian, it may be regarded as a right $\mathbb{Z}(M_n/M_n')$-module in the usual way, where the module action comes from conjugation in $M_n$. Write $A_n = F_n/F_n'$ for the free abelian group of rank $n$ and $a_i = f_iF_n'$, $i = 1, \ldots, n$. Thus, the set $\{a_1, \ldots, a_n\}$ is a free generating set for $A_n$. The natural epimorphism $\pi : M_n \to A_n$ induces an isomorphism from $M_n/M_n'$ to $A_n$. So, we may regard $M_n'$ as a right $\mathbb{Z}A_n$-module. For $w \in M_n'$ and $s \in \mathbb{Z}A_n$, we write $w^s$ to denote the image of $w$ under the action of $s$. This notation is consistent with our notation for conjugation. A set of generators of $M_n'$ is given by $[x_i, x_j]^{a_1^{\kappa_1} \cdots a_n^{\kappa_n}}$, where $i < j$ and $\kappa_\mu \in \mathbb{Z}$, $\mu \in \{1, \ldots, n\}$. Since $M_n'$ is abelian, any element of $M_n'$ may be written in the form $\prod_{1 \leq i < j \leq n} [x_i, x_j]^{P_{ij}}$, where $P_{ij} \in \mathbb{Z}A_n$. Since $M_n$ is metabelian, it satisfies the identity $[u, v, w][v, w, u][w, u, v] = 1$ for all $u, v, w \in M_n$. The commutators in this paper are taken as left-normed. In the case where $u = x_i$, $v = x_j$ and $w = x_k$, we have $[x_i, x_j]^{(a_k - 1)}[x_j, x_k]^{(a_i - 1)}[x_k, x_i]^{(a_j - 1)} = 1$. Thus, for $n \geq 3$, the above mentioned generating set for $M_n'$ is not a free generating set.

### 4.3.1 Modules. Embedding.

The natural epimorphism from $M_n$ onto $M_n/M_n'$ induces a group homomorphism $\tilde{\sigma}_n$ from $\mathrm{Aut}(M_n)$ into $\mathrm{GL}_n(\mathbb{Z})$. Since the natural homomorphism from $\mathrm{Aut}(F_n)$ into $\mathrm{GL}_n(\mathbb{Z})$ is surjective [1, Theorem 3.1] and $M_n/M_n' \cong F_n/F_n'$, we have $\tilde{\sigma}_n$ is surjective and the kernel of $\tilde{\sigma}_n$ is equal to $\mathrm{IA}(M_n)$. Then, $\mathrm{Aut}(M_n) = T_n\,\mathrm{IA}(M_n)$. Since $\mathrm{IA}(M_n)$ is normal, $\mathrm{Aut}(M_n)$ acts by conjugation on $\mathrm{IA}(M_n)$. Since $\mathrm{IA}(M_n)$ acts trivially by conjugation on each $\mathcal{L}^c(\mathrm{IA}(M_n))$, with $c \geq 2$, we have each $\mathcal{L}^c(\mathrm{IA}(M_n))$ is a $\mathbb{Z}T_n$-module. Let $\overline{\phi} = \phi\ \mathrm{I}_{c+1}\mathrm{A}(M_n)$, with $\phi \in \mathrm{I}_c\mathrm{A}(M_n)$, and let $g \in \mathrm{GL}_n(\mathbb{Z})$. There exists $t_g \in T_n$ such that $\tilde{\sigma}_n(t_g) = g$. Let $t'_g \in T_n$ such that $\tilde{\sigma}_n(t'_g) = g$. Thus, there exists $\psi \in \mathrm{IA}(M_n)$ such that $t_g = t'_g\psi$. Since $[\mathrm{I}_c\mathrm{A}(M_n), \mathrm{IA}(M_n)] \subseteq \mathrm{I}_{c+1}\mathrm{A}(M_n)$ and $\mathrm{I}_{c+1}\mathrm{A}(M_n)$ is normal in $\mathrm{Aut}(M_n)$, we have $(t_g\phi t'_g)\mathrm{I}_{c+1}\mathrm{A}(M_n) = (t'_g\phi(t'_g)^{-1})\mathrm{I}_{c+1}\mathrm{A}(M_n)$. Hence, we may define an action of $\mathrm{GL}_n(\mathbb{Z})$ on $\mathcal{L}^c(\mathrm{IA}(M_n))$ by $g * \overline{\phi} = \overline{t_g\phi t_g^{-1}}$ for all $g \in \mathrm{GL}_n(\mathbb{Z})$, $\overline{\phi} = \phi\ \mathrm{I}_{c+1}\mathrm{A}(M_n)$, with $\phi \in \mathrm{I}_c\mathrm{A}(M_n)$, and $\tilde{\sigma}_n(t_g) = g$. It turns out that the above action is a left action of $\mathrm{GL}_n(\mathbb{Z})$ on $\mathcal{L}^c(\mathrm{IA}(M_n))$. The action of $\mathrm{GL}_n(\mathbb{Z})$ commutes with the multiplication of the elements of $\mathbb{Z}$. So, $\mathcal{L}^c(\mathrm{IA}(M_n))$ is a $\mathbb{Z}\mathrm{GL}_n(\mathbb{Z})$-module for all $c \geq 2$. In this paper, all $\mathbb{Z}\mathrm{GL}_n(\mathbb{Z})$-modules are regarded as left $\mathbb{Z}\mathrm{GL}_n(\mathbb{Z})$-modules.

We write $\mathrm{gr}_{1,\mathbb{Q}}(M_n)$ for the tensor product of $\mathbb{Q}$ with $\mathrm{gr}_1(M_n)(\cong \mathrm{gr}_1(F_n))$ over $\mathbb{Z}$. Since $\mathrm{gr}_1(M_n)$ is a free $\mathbb{Z}$-module with a free generating set $\{\overline{x}_1, \ldots, \overline{x}_n\}$, where $\overline{x}_i = x_iM'_n$ for $i \in [n] = \{1, \ldots, n\}$, we may regard $\mathrm{gr}_1(M_n) \subseteq \mathrm{gr}_{1,\mathbb{Q}}(M_n)$. Furthermore, the set $\{\overline{x}_1, \ldots, \overline{x}_n\}$ is a $\mathbb{Q}$-basis of $\mathrm{gr}_{1,\mathbb{Q}}(M_n)$. Thus, we may think of the elements of $\mathrm{gr}_{1,\mathbb{Q}}(M_n)$ as a $\mathbb{Q}$-linear combination of $\overline{x}_1, \ldots, \overline{x}_n$. We identify $\mathrm{Aut}(\mathrm{gr}_{1,\mathbb{Q}}(M_n))$ with $\mathrm{GL}_n(\mathbb{Q})$ with respect to the $\mathbb{Q}$-basis $\{\overline{x}_1, \ldots, \overline{x}_n\}$. The group $\mathrm{GL}_n(\mathbb{Q})$ acts naturally on $\mathrm{gr}_{1,\mathbb{Q}}(M_n)$ by $g\overline{x}_j = \sum_{i=1}^n g_{ij}\overline{x}_i$ for $j \in [n]$, where $g = (g_{ij}) \in \mathrm{GL}_n(\mathbb{Q})$ for $i, j \in [n]$. For $c \geq 2$, let $\mathrm{gr}_{c,\mathbb{Q}}(M_n)$ denote the tensor product of $\mathbb{Q}$ with $\mathrm{gr}_c(M_n)$ over $\mathbb{Z}$. As mentioned before, the set $\mathcal{M}_{n,c}$ is a $\mathbb{Q}$-basis of $\mathrm{gr}_{c,\mathbb{Q}}(M_n)$. Since $\mathrm{gr}_c(M_n)$ is a free $\mathbb{Z}$-module with a free generating set $\mathcal{M}_{n,c}$, we may regard $\mathrm{gr}_c(M_n) \subseteq \mathrm{gr}_{c,\mathbb{Q}}(M_n)$ and we may consider the elements of $\mathrm{gr}_{c,\mathbb{Q}}(M_n)$ as a $\mathbb{Q}$-linear combination of elements of $\mathcal{M}_{n,c}$. The action of $\mathrm{GL}_n(\mathbb{Q})$ on $\mathrm{gr}_{1,\mathbb{Q}}(M_n)$ can be extended diagonally on $\mathrm{gr}_{c,\mathbb{Q}}(M_n)$, for $c \geq 2$, subject to

$$g\overline{[x_{j_1}, \ldots, x_{j_c}]} = g([x_{j_1}, \ldots, x_{j_c}] + \gamma_{c+1}(M_n)) = [g\overline{x}_{j_1}, \ldots, g\overline{x}_{j_c}] + \gamma_{c+1}(M_n) = \overline{[g\overline{x}_{j_1}, \ldots, g\overline{x}_{j_c}]}$$

with $j_1 > j_2 \leq j_3 \leq \cdots \leq j_c$. It turns out that the above action is a left action of $\mathrm{GL}_n(\mathbb{Q})$ on $\mathrm{gr}_{c,\mathbb{Q}}(M_n)$. The action of $\mathrm{GL}_n(\mathbb{Q})$ commutes with the multiplication of the elements of $\mathbb{Q}$. So, each $\mathrm{gr}_{c,\mathbb{Q}}(M_n)$ is a $\mathbb{Q}\mathrm{GL}_n(\mathbb{Q})$-module.

Let $L_n$ be a free Lie algebra (over $\mathbb{Q}$) of rank $n$, with $n \geq 2$, with a free generating set $\{z_1, \ldots, z_n\}$ and $L'_n$ be the derived algebra of $L_n$. Write $M_{n,L} = L_n/L''_n$, where $L''_n = (L'_n)'$, and $z'_j = z_j + L''_n$, $j \in [n]$. Thus, $M_{n,L}$ is a free metabelian Lie algebra of rank $n$. For $c \in \mathbb{N}$, let $M^c_{n,L}$ be the vector subspace of $M_{n,L}$ spanned by all Lie commutators $[z'_{j_1}, \ldots, z'_{j_c}]$, with $j_1, \ldots, j_c \in [n]$. So, $M_{n,L} = \bigoplus_{c\geq 1} M^c_{n,L}$. We point out that each $M^c_{n,L}$ becomes a $\mathbb{Q}\mathrm{GL}_n(\mathbb{Q})$-module in a natural way. Each $M^c_{n,L}$ is isomorphic to $\mathrm{gr}_{c,\mathbb{Q}}(M_n)$ as a $\mathbb{Q}\mathrm{GL}_n(\mathbb{Q})$-module by an isomorphism sending $[z'_{j_1}, \ldots, z'_{j_c}]$ to $[\overline{x}_{j_1}, \ldots, \overline{x}_{j_c}] + \gamma_{c+1}(M_n)$ for all $j_1, \ldots, j_c \in [n]$. The group $\mathrm{GL}_n(\mathbb{Q})$ acts on $\mathrm{gr}_{c,\mathbb{Q}}(M_n)^{\oplus n}$ by $g \bullet (\overline{u}_1, \ldots, \overline{u}_n) = (g\overline{u}_1, \ldots, g\overline{u}_n)g^{-1}$, where $g \in \mathrm{GL}_n(\mathbb{Q})$ and $\overline{u}_i \in \mathrm{gr}_{c,\mathbb{Q}}(M_n)$ for $i \in [n]$. Here $g\overline{u}_i$ means the canonical action of $g$ on $\mathrm{gr}_{c,\mathbb{Q}}(M_n)$ and the multiplication of a $1 \times n$ and an $n \times n$ matrix. Clearly, each $\mathrm{gr}_{c,\mathbb{Q}}(M_n)^{\oplus n}$ is a $\mathbb{Q}\mathrm{GL}_n(\mathbb{Q})$-module.

For $c \geq 2$, let $\chi_c$ be the mapping from $\mathrm{I}_c\mathrm{A}(M_n)$ to $\mathrm{gr}_{c,\mathbb{Q}}(M_n)^{\oplus n}$ defined by $\chi_c(\phi) = (\overline{u}_1, \ldots, \overline{u}_n)$ for all $\phi \in \mathrm{I}_c\mathrm{A}(M_n)$, with $\phi(x_i) \equiv x_iu_i \pmod{\gamma_{c+1}(\mathrm{M_n})}$, $u_i \in \gamma_c(M_n)$, $i \in [n]$.

It is easily verified that $\chi_c$ is a group homomorphism and the kernel of $\chi_c$ is equal into $\mathrm{I}_{c+1}\mathrm{A}(M_n)$. Therefore, via the monomorphism $\overline{\chi}_c$ induced by $\chi_c$, $\mathcal{L}^c(\mathrm{IA}(M_n))$ is isomorphic to a subgroup of $\mathrm{gr}_{c,\mathbb{Q}}(M_n)^{\oplus n}$ in a natural way. It is easily checked that $\overline{\chi}_c(g * \overline{\phi}) = g \bullet (\overline{u}_1, \ldots, \overline{u}_n)$ for all $g \in \mathrm{GL}_n(\mathbb{Z})$ and $\phi \in \mathrm{I}_c\mathrm{A}(M_n)$, with $\phi(x_i) \equiv x_i u_i \gamma_{c+1}(M_n)$, $u_i \in \gamma_c(M_n)$, $i \in [n]$.

### 4.3.2 Fox derivatives. Criterion.

For $j \in [n]$, the (left) Fox derivative associated with $f_j$ is the linear map $D_j : \mathbb{Z}F_n \to \mathbb{Z}F_n$ satisfying the conditions: $D_j(f_j) = 1$, $D_j(f_i) = 0$ for $i \neq j$ and $D_j(uv) = D_j(u) + uD_j(v)$ for all $u, v \in F_n$. Since $D_j(1) = 0$, we have $D_j(u^{-1}) = -u^{-1}D_j(u)$ for all $u \in F_n$. It is easily verified that $D_j([u,v]) = u^{-1}(v^{-1}-1)D_j(u) + u^{-1}v^{-1}(u-1)D_j(v)$ for all $u, v \in F_n$. Let $\mathcal{A} = \mathrm{Ker}(\varepsilon_F : \mathbb{Z}F_n \to \mathbb{Z})$ be the augmentation ideal of $\mathbb{Z}F_n$. The kernel of $\varepsilon_F$ is a free left $\mathbb{Z}F_n$-module with basis $\{f_i - 1 : i \in [n]\}$ (see, for example, [10]). By an inductive argument on $c$, with $c \geq 2$, we have $[u,v]-1 \in \mathcal{A}^c$ for all $u \in \gamma_{c-1}(F_n)$ and $v \in F_n$, and $D_j(w) \in \mathcal{A}^{c-1}$ for all $w \in \gamma_c(F_n)$ and $j \in [n]$. The natural epimorphisms $\pi : M_n \to A_n$, $\pi_A : F_n \to A_n$ and $\pi_M : F_n \to M_n$ extend to ring epimorphisms $\pi : \mathbb{Z}M_n \to \mathbb{Z}A_n$, $\pi_A : \mathbb{Z}F_n \to \mathbb{Z}A_n$ and $\pi_M : \mathbb{Z}F_n \to \mathbb{Z}M_n$. The kernels of $\pi_A$ and $\pi_M$ are the ideals of $\mathbb{Z}F_n$ generated by the elements $u - 1$, with $u \in F_n'$ and $u \in F_n''$, respectively. It follows that $D_j(w) \in \mathrm{Ker}(\pi_A)$ for all $w \in F_n''$ and so, $D_j(\mathrm{Ker}\pi_M) \subseteq \mathrm{Ker}\pi_A$. For $j \in [n]$, we define $\partial_j : \mathbb{Z}M_n \to \mathbb{Z}A_n$ by $\partial_j(\pi_M(u)) = \pi_A(D_j(u))$ for all $u \in \mathbb{Z}F_n$. Since $D_j(\mathrm{Ker}\pi_M) \subseteq \mathrm{Ker}\pi_A$, we get $\partial_j$ is well defined for all $j \in [n]$. Thus, $\partial_j$ is a linear map for all $j \in [n]$. From the definitions, we have $\partial_j(x_j) = 1$, $\partial_j(x_i) = 0$ for $i \neq j$ and $\partial_j(uv) = \partial_j(u) + \pi_A(u)\partial_j(v)$ for all $u, v \in M_n$. Moreover,

$$\partial_j([u,v]) = \pi_A(u^{-1})(\pi_A(v^{-1})-1)\partial_j(u) + \pi_A(u^{-1})\pi_A(v^{-1})(\pi_A(u)-1)\partial_j(v) \tag{3.3}$$

for all $u, v \in M_n$ and $j \in [n]$. In particular, by the equation (3.3), for $j \neq k$, $\partial_j([x_j, x_k]) = a_j^{-1}a_k^{-1}(1-a_k)$ and $\partial_k([x_j, x_k] = a_j^{-1}a_k^{-1}(a_j - 1)$. Note that $\sum_{j=1}^n \partial_j(x_i w_i)(a_j - 1) = a_i - 1$ for $i \in [n]$ and $w \in M_n'$. For any IA-endomorphism $\phi$ of $M_n$, we define the Jacobian matrix $J(\phi) = (a_{ij}) \in M_{n\times n}(\mathbb{Z}A_n)$, where $a_{ij} = \partial_j(\phi(x_i))$ for all $i, j \in [n]$. It is well known that $\phi \in \mathrm{IA}(M_n)$ if and only if the determinant $\det(J(\phi))$ of $J(\phi)$ is $a_1^{i_1} \cdots a_n^{i_n}$ (see, [3, Proposition 2], [2, Theorem 1], [13], [10, Theorem 1.16]). For $s \in \mathbb{Z}A_n$, write $s = \sum_i m_i s_i$, where $m_i \in \mathbb{Z}$ and $s_i \in A_n$ for each $i$, and define $s^* = \sum_i m_i s_i^{-1}$. Note that $(a_{i_1} \cdots a_{i_\nu})^* = a_{i_1}^{-1} \cdots a_{i_\nu}^{-1}$ for all $i_1, \ldots, i_\mu \in [n]$. The map $* : \mathbb{Z}A_n \to \mathbb{Z}A_n$ sending $s$ to $s^*$ is an involutory linear map. For $w \in M_n'$ and $s \in \mathbb{Z}A_n$, it is easily verified that $\partial_j(w^s) = s^*\partial_j(w)$ for all $j \in [n]$.

### 4.3.3 A generating set for $\mathcal{L}^c(\mathrm{IA}(M_n))$

Let $\Sigma = \mathrm{Ker}(\varepsilon_A : \mathbb{Z}A_n \to \mathbb{Z})$ be the augmentation ideal of $\mathbb{Z}A_n$. Since $\pi_A(\mathcal{A}^m) = \Sigma^m$ for all $m \geq 1$, we get $\partial_j(w') \in \Sigma^{c-1}$ for all $w' \in \gamma_c(M_n)$ and $j \in [n]$. As observed in [3, Lemma 6], if $\phi \in \mathrm{IA}(M_n)$, then $\phi \in \mathcal{L}^c(\mathrm{IA}(M_n))$, with $c \geq 3$, if and only if $\det(J(\phi)) = 1$. Let $\phi \in \mathrm{I}_c\mathrm{A}(M_n)$, with $c \geq 3$, and so, $\phi(x_i) = x_i u_i$ with $u_i \in \gamma_c(M_n)$ and $i \in [n]$. Write each $u_i = \prod_{1\leq j<k\leq n}[x_j,x_k]^{P_{ijk}}v_i$, where $v_i \in \gamma_{c+1}(M_n)$, $i \in [n]$ and each $P_{ijk}$ is a monomial in the $a_1-1, a_2-1, \ldots, a_n-1$ of total degree $c-2$. Thus, $J(\phi) \equiv I_n \pmod{\Sigma^{c-1}}$. For distinct $i,j,k \in [n]$ and non-negative integers $r_1,\ldots,r_n$, with $r_1+\cdots+r_n = c-2$, let $\tau_{ijk,(r_1,\ldots,r_n)}$ be the IA-endomorphism of $M_n$ satisfying the conditions $\tau_{ijk,(r_1,\ldots,r_n)}(x_i) = x_i[x_j,x_k]^{(a_1-1)^{r_1}\cdots(a_n-1)^{r_n}}$ and $\tau_{ijk,(r_1,\ldots,r_n)}(x_r) = x_r$, $r \neq i$. Since $\det(J(\tau_{ijk,(r_1,\ldots,r_n)})) = 1$, we have $\tau_{ijk,(r_1,\ldots,r_n)} \in \mathrm{I}_c\mathrm{A}(M_n)$. For $j_1,\ldots,j_{c-1} \in [n]$, with $j_1 \neq j_2$, we write $u(j_1,\ldots,j_{c-1}) = [x_{j_1},\ldots,x_{j_{c-1}}]$, and let $\xi_{u(j_1,\ldots,j_{c-1})}$ be the inner automorphism of $M_n$ on $u(j_1,\ldots,j_{c-1})$. That is, for any $j \in [n]$, $\xi_{u(j_1,\ldots,j_{c-1})}(x_j) = x_j[x_j,u(j_1,\ldots,j_{c-1})]$. Clearly, $\xi_{u(j_1,\ldots,j_{c-1})} \in \mathrm{I}_c\mathrm{A}(M_n)$.

For $c \geq 3$, let $P,Q \in \mathbb{Z}A_n$ be monomials in the $a_1-1,\ldots,a_n-1$ of total degree $c-2$ and $c-3$, respectively. For distinct $i,j,k \in [n]$, let $B_{ikj}(P)$, and $B_{ij}(Q)$ be the IA-endomorphisms of $M_n$ satisfying the conditions

$$\begin{aligned} B_{ikj}(P)(x_i) &= x_i[x_i,x_j]^{a_i^{-1}a_j^{-1}a_k^{-1}P}[x_k,x_j]^{a_k^{-2}a_j^{-1}P}, \\ B_{ikj}(P)(x_k) &= x_k[x_i,x_j]^{-a_i^{-2}a_j^{-1}P}[x_k,x_j]^{-a_i^{-1}a_j^{-1}a_k^{-1}P} \text{ and} \\ B_{ikj}(P)(x_r) &= x_r, \text{with } r \neq i,k, \end{aligned}$$

and

$$\begin{aligned} B_{ij}(Q)(x_i) &= x_i[x_i,x_j]^{-a_i^{-2}a_j^{-2}(a_i-1)Q}, \\ B_{ij}(Q)(x_j) &= x_j[x_i,x_j]^{-a_i^{-2}a_j^{-2}(a_j-1)Q} \text{ and} \\ B_{ij}(Q)(x_s) &= x_s, \text{ with } s \neq i,j. \end{aligned}$$

Since $\det(J(B_{ikj}(P))) = \det(J(B_{ij}(Q))) = 1$, we have $B_{ikj}(P), B_{ij}(Q) \in \mathrm{I}_c\mathrm{A}(M_n)$. Thus, $B_{ikj}(P)(x_i) = x_i[x_i,x_j]^P[x_k,x_j]^P v_i$, $B_{ikj}(P)(x_k) = x_k[x_j,x_i]^P[x_j,x_k]^P v_k$, $B_{ikj}(P)(x_r) = x_r$ with $r \neq i,k$ and $v_i, v_k \in \gamma_{c+1}(M_n)$. Furthermore, $B_{ij}(Q)(x_i) = x_i[x_i,[x_i,x_j]]^Q v_i'$, $B_{ij}(Q)(x_j) = x_j[x_j,[x_i,x_j]]^Q v_k'$ and $B_{ij}(Q)(x_s) = x_s$, with $s \neq i,j$ and $v_i', v_k' \in \gamma_{c+1}(M_n)$.

**Lemma 4** *For $c \geq 3$, $\mathcal{L}^c(\mathrm{IA}(M_n))$ is generated as a $\mathbb{Z}\mathrm{GL}_n(\mathbb{Z})$-module by all $\overline{\tau}_{123,(r_1,\ldots,r_n)}$, $\overline{B}_{123}(P)$ and $\overline{B}_{12}(Q)$ with $r_1+\cdots+r_n = c-2$, $P$ and $Q$ are monomials in $a_1-1,\ldots,a_n-1$ of total degree $c-2$ and $c-3$, respectively.*

*Proof.* It follows from [3, Theorem 3 and Proof of Lemma 7] that, for $c \geq 3$, the $\mathbb{Z}$-module $\mathcal{L}^c(\mathrm{IA}(M_n))$ is generated by all $\overline{\tau}_{ijk,(r_1,\ldots,r_n)}$, $\overline{B}_{ikj}(P)$ and $\overline{B}_{ij}(Q)$. For suitable choices of $g_1,g_2,g_3 \in \mathrm{GL}_n(\mathbb{Z})$, we have $g_1 * \overline{\tau}_{123,(s_1,\ldots,s_n)} = \overline{\tau}_{ijk,(r_1,\ldots,r_n)}$, $g_2 * \overline{B}_{123}(P') = \overline{B}_{ikj}(P)$ and $g_3 * \overline{B}_{12}(Q') = \overline{B}_{ij}(Q)$ and so, we obtain the required result. □

For the proof of Theorem 1 (3), the following result gives us a suitable generating set for $\mathcal{L}^c(\mathrm{IA}(M_n))$ as a $\mathbb{Z}\mathrm{GL}_n(\mathbb{Z})$-module.

**Proposition 2** *Let $M_n$ be a free metabelian group of finite rank $n$, with $n \geq 4$. Then, for any $c \geq 4$, the $\mathbb{Z}$-module $\mathcal{L}^c(\mathrm{IA}(M_n))$ is generated by $\overline{\tau}_{123,(r_1,\ldots,r_n)}$ and $\overline{\xi}_{v(k_1,\ldots,k_{c-1})}$, with*

$v(j_1,\ldots,j_{c-1}) = [x_{k_1},\ldots,x_{k_{c-1}}]$, $r_1+\cdots+r_n = c-2$, $k_1,\ldots,k_{c-1}\in[n]$ *and* $k_1 > k_2 \leq k_3 \leq \cdots \leq k_{c-1}$ *as a* $\mathbb{Z}\mathrm{GL}_n(\mathbb{Z})$*-module.*

*Proof.* Let $\Omega_{n,c}$ be the $\mathbb{Z}\mathrm{GL}_n(\mathbb{Z})$-submodule of $\mathcal{L}^c(\mathrm{IA}(M_n))$ generated by all $\overline{\tau}_{123,(r_1,\ldots,r_n)}$, with $r_1,\ldots,r_n\geq 0$, $r_1+\cdots+r_n = c-2$ and $\overline{\xi}_{v(k_1,\ldots,k_{c-1})}$, with $v(k_1,\ldots,k_{c-1}) = [x_{k_1},\ldots,x_{k_{c-1}}]$, $k_1,\ldots,k_{c-1}\in[n]$ and $k_1>k_2\leq k_3\leq\cdots\leq k_{c-1}$. Since $\mathcal{M}_{n,c-1}$ is a $\mathbb{Z}$-basis of $\mathrm{gr}_{c-1}(M_n)$, we get, modulo $\gamma_c(M_n)$, each $u(j_1,\ldots,j_{c-1})$, with $j_1,\ldots,j_{c-1}\in[n]$, $j_1\neq j_2$, is written as a product of $v(k_1,\ldots,k_{c-1})$, with $k_1,\ldots,k_{c-1}\in[n]$ and $k_1>k_2\leq k_3\leq\cdots\leq k_{c-1}$. Therefore, each $\overline{\xi}_{u(j_1,\ldots,j_{c-1})}$ is a product of $\overline{\xi}_{u(k_1,\ldots,k_{c-1})}$ and so, $\overline{\xi}_{u(j_1,\ldots,j_{c-1})}\in\Omega_{n,c}$ for all $j_1,\ldots,j_{c-1}\in[n]$, with $j_1\neq j_2$. By Lemma 4, to prove that $\mathcal{L}^c(\mathrm{IA}(M_n)) = \Omega_{n,c}$ for all $c\geq 3$, it is enough to show that $\overline{B}_{123}(P)$ and $\overline{B}_{12}(Q)$ belong to $\Omega_{n,c}$ for all monomials $P$ and $Q$ of total degree $c-2$ and $c-3$, respectively.

Throughout the proof, for $c\geq 3$ and non-negative integers $r_1,\ldots,r_n$, $s_1,\ldots,s_n$, with $r_1+\cdots r_n = c-2$ and $s_1+\cdots+s_n = c-3$, we write $B_{123}(r_1,\ldots,r_n)$ and $B_{12}(s_1,\ldots,s_n)$ instead of $B_{123}(P)$ and $B_{12}(Q)$, respectively, where $P = (a_1-1)^{r_1}\cdots(a_n-1)^{r_n}$ and $Q = (a_1-1)^{s_1}\cdots(a_n-1)^{s_n}$. For non-negative integers $s_1,\ldots,s_n$, with $s_1+\cdots+s_n = c-3$, let $u(1,2,s_1,\ldots,s_n) = [x_1,x_2,\ {}_{s_1}x_1,\ldots,\ {}_{s_n}x_n]$. Let $\xi_{u(1,2,s_1,\ldots,s_n)}$ be the inner automorphism of $M_n$ on $u(1,2,s_1,\ldots,s_n)$. Thus, as shown above, $\overline{\xi}_{u(1,2,s_1,\ldots,s_n)}\in\Omega_{n,c}$. Since $M_n$ is metabelian, we have

$$[x_1,[x_1,x_2],\ {}_{s_1}x_1,\ldots,\ {}_{s_n}x_n] = [x_1,[x_1,x_2,\ {}_{s_1}x_1,\ldots,\ {}_{s_n}x_n]] = [x_1,u(1,2,s_1,\ldots,s_n)].$$

By direct calculations,

$$\overline{\xi}_{u(1,2,s_1,\ldots,s_n)}\Big(\prod_{\lambda=3}^{n}(\overline{\tau}_{\lambda 12,(s_1,\ldots,s_{\lambda-1},s_\lambda+1,s_{\lambda+1},\ldots,s_n)})^{-1}\Big) = \overline{B}_{12}(s_1,\ldots,s_n). \tag{3.4}$$

By the proof of Lemma 4, we have $\overline{\tau}_{\lambda 12,(s_1,\ldots,s_{\lambda-1},s_\lambda+1,s_{\lambda+1},\ldots,s_n)}\in\Omega_{n,c}$ for all $\lambda\in\{3,\ldots,n\}$. Since $\overline{\xi}_{u(1,2,s_1,\ldots,s_n)}\in\Omega_{n,c}$, by the equation (3.4), we get $\overline{B}_{12}(s_1,\ldots,s_n)\in\Omega_{n,c}$ for all non-negative integers $s_1,\ldots,s_n$, with $s_1+\cdots+s_n = c-3$.

For non-negative integers $r_1,\ldots,r_n$, with $r_1+\cdots+r_n = c-2$, we define $\delta_{123,(r_1,\ldots,r_n)} = B_{123}(r_1,\ldots,r_n)\tau_{213,(r_1,\ldots,r_n)}\tau^{-1}_{123,(r_1,\ldots,r_n)}$. By working modulo $\gamma_{c+1}(M_n)$, $\delta_{123,(r_1,\ldots,r_n)}$ satisfies the conditions

$$\begin{array}{rcl}\delta_{123,(r_1,\ldots,r_n)}(x_1) & = & x_1[x_1,x_3,\ {}_{r_1}x_1,\ldots,\ {}_{r_n}x_n],\\ \delta_{123,(r_1,\ldots,r_n)}(x_2) & = & x_2[x_3,x_2,\ {}_{r_1}x_1,\ldots,\ {}_{r_n}x_n] \text{ and}\\ \delta_{123,(r_1,\ldots,r_n)}(x_t) & = & x_t \text{ for } t\neq 1,2.\end{array}$$

Let $\beta$ be the automorphism of $M_n$ satisfying the conditions $\beta(x_2) = x_2x_1$ and $\beta(x_\mu) = x_\mu$, $\mu\neq 2$. By direct calculations, for all $r_1,r_3,\ldots,r_n\geq 0$, with $r_1+r_3+\cdots+r_n = c-2$, we have

$$\overline{\tau}_{213,(r_1,0,r_3,\ldots,r_n)}(\overline{\tau}_{123,(r_1,0,r_3,\ldots,r_n)})^{-1}(\beta^{-1}*\overline{\tau}_{123,(r_1,0,r_3,\ldots,r_n)}) = \overline{\delta}_{123,(r_1,0,r_3,\ldots,r_n)} \tag{3.5}$$

and so, $\overline{\delta}_{123,(r_1,0,r_3,\ldots,r_n)}\in\Omega_{n,c}$. Thus, we may assume that $r_2\geq 1$. First, we consider the case $r_2 = 1$. For non-negative integers $r_1,r_3,\ldots,r_n$ and $r_1+r_3+\cdots+r_n = c-3$, we write $\phi_{1,1} = \tau^{-1}_{123,(r_1,1,r_3,\ldots,r_n)}\tau^{-1}_{123,(r_1+1,0,r_3,\ldots,r_n)}$ and $\phi_{2,1} = \tau_{213,(r_1,1,r_3,\ldots,r_n)}\tau_{213,(r_1+1,0,r_3,\ldots,r_n)}$. By direct calculations, we have

$$(\overline{\delta}_{123,(r_1+1,0,r_3,\ldots,r_n)})^{-1}\overline{\phi}_{2,1}\overline{\phi}_{1,1}(\beta^{-1} * \overline{\tau}_{123,(r_1,1,r_3,\ldots,r_n)}) = \overline{\delta}_{123,(r_1,1,r_3,\ldots,r_n)} \tag{3.6}$$

for all non-negative integers $r_1, r_3, \ldots, r_n$, with $r_1 + r_3 + \cdots + r_n = c - 3$. By the equation (3.6), we have $\overline{\delta}_{123,(r_1,1,r_3,\ldots,r_n)} \in \Omega_{n,c}$. Next, we assume that $r_2 \geq 2$. Define

$B = \prod_{k=1}^{r_2-1} B_{12,(r_1+k-1,r_2-k-1,r_3+1,r_4,\ldots,r_n)}$, $\psi_{1,1} = \prod_{k=1}^{r_2-1}(\tau_{123,(r_1+k+1,r_2-k-1,r_3,\ldots,r_n)})^{-1}$, $\psi_{2,2} = \prod_{k=1}^{r_2-1} \tau_{213,(r_1+k-1,r_2-k+1,r_3,\ldots,r_n)}$, $\psi_{1,3} = \prod_{k=0}^{r_2}(\tau_{123,(r_1+k,r_2-k,r_3,\ldots,r_n)})^{-1}$ and $\psi_{2,4} = \prod_{k=0}^{r_2} \tau_{213,(r_1+k,r_2-k,r_3,\ldots,r_n)}$.

It is easily verified that, for all $r_1, \ldots, r_n \geq 0$, with $r_1 + \cdots + r_n = c - 2$ and $r_2 \geq 2$,

$$\overline{B}\ \overline{\psi}_{2,2}\overline{\psi}_{1,1}(\overline{\delta}_{123,(r_1+r_2,0,r_3,\ldots,r_n)})^{-1}\overline{\psi}_{2,4}\overline{\psi}_{1,3}(\beta^{-1} * \overline{\tau}_{123,(r_1,r_2,\ldots,r_n)}) = \overline{\delta}_{123,(r_1,\ldots,r_n)}. \tag{3.7}$$

By the equation (3.7), we have $\overline{\delta}_{123,(r_1,1,r_3,\ldots,r_n)} \in \Omega_{n,c}$. Hence, by the equations (3.4) – (3.7), we obtain $\overline{\delta}_{123,(r_1,\ldots,r_n)} \in \Omega_{n,c}$ for all $r_1, \ldots, r_n \geq 0$, with $r_1 + \cdots + r_n = c - 2$. Since $\delta_{123,(r_1,\ldots,r_n)} = B_{123}(r_1, \ldots, r_n)\tau_{213,(r_1,\ldots,r_n)}\tau^{-1}_{123,(r_1,\ldots,r_n)}$, we get $\overline{B}_{123}(r_1, \ldots, r_n) \in \Omega_{n,c}$ for all $r_1, \ldots, r_n \geq 0$ and $r_1 + r_2 + \cdots + r_n = c - 2$. □

Let $n \geq 4$ and $c \geq 3$. For $i \in [n]$ and $i_1, \ldots, i_c \in [n] \setminus \{i\}$, with $i_1 \neq i_2$, let $\tau_{i,(i_1,\ldots,i_c)}$ be the IA-endomorphism of $M_n$ satisfying the conditions $\tau_{i,(i_1,\ldots,i_c)}(x_i) = x_i[x_{i_1}, \ldots, x_{i_c}]$ and $\tau_{i,(i_1,\ldots,i_c)}(x_r) = x_r$, $r \neq i$. It is easily verified that $\tau_{i,(i_1,\ldots,i_c)}$ is an automorphism of $M_n$. Furthermore, for $i \in [n]$ and $i_1, \ldots, i_c \in [n] \setminus \{i\}$, with $i_1 \neq i_2$, we write $\overline{u}_i(i_1, \ldots, i_c) = [\overline{x}_{i_1}, \ldots, \overline{x}_{i_c}] + \gamma_{c+1}(M_n)$. Then, $\overline{\chi}_c(\overline{\tau}_{i,(i_1,\ldots,i_c)})$ is the $n$-tuple with $\overline{u}_i(i_1, \ldots, i_c)$ in the $i$-th position and the rest are zero. Let $P_c$ be the $\mathbb{Q}\mathrm{GL}_n(\mathbb{Q})$-submodule of $\mathrm{gr}_{c,\mathbb{Q}}(M_n)^{\oplus n}$ generated by the set $\mathcal{P}_c = \{\overline{\chi}_c(\overline{\tau}_{1,(i_1,\ldots,i_c)}) : i_1 \neq i_2; i_1, \ldots, i_c \in [n] \setminus \{1\}\}$. Thus, each $\overline{\chi}_c(\overline{\tau}_{i,(i_1,\ldots,i_c)}) \in P_c$.

For $u \in M_n'$, write $\xi_u \in \mathrm{IA}(M_n)$ which satisfies the conditions $\xi_u(x_i) = x_i[x_i, u]$, $i \in [n]$. Note that $\xi_{u^{-1}} = (\xi_u)^{-1}$ for all $u \in M_n'$. Let $\mathcal{Q}_c = \{\overline{\chi}_c(\overline{\xi}_u) = (\overline{[x_1, u]}, \ldots, \overline{[x_n, u]}) : u \in \gamma_{c-1}(M_n)\}$. It is easily verified that if $c \geq 3$, $\psi \in \mathrm{Aut}(M_n)$, $u, u_1, u_2 \in \gamma_{c-1}(M_n)$ and $a \in \mathbb{Z}$, then $\psi\xi_u\psi^{-1} = \xi_{\psi(u)}$, $\xi_{u_1}\xi_{u_2} = \xi_{u_1u_2}$ and $a\overline{\xi}_u = \overline{\xi}_{u^a}$. Since $\overline{\chi}_c$ is a $\mathbb{Z}\mathrm{GL}_n(\mathbb{Z})$-monomorphism, we have $\mathcal{Q}_c$ is a $\mathbb{Z}\mathrm{GL}_n(\mathbb{Z})$-module. Let $Q_c$ be the $\mathbb{Q}\mathrm{GL}_n(\mathbb{Q})$-submodule of $\mathrm{gr}_{c,\mathbb{Q}}(M_n)^{\oplus n}$ generated by the set $\mathcal{Q}_c$.

Let $\eta$ be the IA-endomorphism of $M_n$ satisfying the conditions $\eta(x_j) = x_j[x_j,\ _{(c-1)}x_1]$ for all $j \in [n]$. Since $\det(J(\eta)) \neq a_1^{k_1} \cdots a_n^{k_n}$ for all $k_1, \ldots, k_n \in \mathbb{Z}$, we have $\eta \notin \mathrm{I_cA}(M_n)$. Let $R_c$ be the $\mathbb{Q}\mathrm{GL}_n(\mathbb{Q})$-submodule of $\mathrm{gr}_{c,\mathbb{Q}}(M_n)^{\oplus n}$ generated by $(0, \overline{[x_2,\ _{(c-1)}x_1]}, \ldots, \overline{[x_n,\ _{(c-1)}x_1]})$.

It follows from [6, Section 2, Proposition 3.5 (i)] that $\mathrm{gr}_{c,\mathbb{Q}}(M_n)^{\oplus n} = P_c \oplus Q_c \oplus R_c$ and $P_c, Q_c$ and $R_c$ are irreducible (rational) $\mathbb{Q}\mathrm{GL}_n(\mathbb{Q})$-modules. Furthermore, $\dim(P_c \oplus Q_c) = \mathrm{rank}(\mathcal{L}^c(\mathrm{IA}(M_n))) = n(c-1)\binom{n+c-2}{n-2} - \binom{n+c-2}{n-1}$ for all $c \geq 3$. By Proposition 2, we obtain the following result.

**Corollary 2** *With the previous notations, for all $n \geq 4$ and $c \geq 3$, $\mathbb{Q} \otimes \overline{\chi}_c(\mathcal{L}^c(\mathrm{IA}(M_n))) = P_c \oplus Q_c$ as vector spaces over $\mathbb{Q}$, where $P_c$ and $Q_c$ are certain irreducible (rational) $\mathbb{Q}\mathrm{GL}_n(\mathbb{Q})$-modules.*

We next need the following technical result.

**Lemma 5** *Let $M_n$ be a free metabelian group of rank $n$, with $n \geq 4$ and let $c \geq 2$. Then,*

1. *For $i_1, \ldots, i_c \in [n] \setminus \{1\}$, $\tau_{1,(i_1,\ldots,i_c)} \in \gamma_{c-1}(\mathrm{IA}(M_n))$.*

2. *For non-negative integers $r_1, r_2, \ldots, r_n$, with $r_1 + r_2 + \cdots + r_n = c-2$, $\tau_{123,(r_1,r_2,\ldots,r_n)} \in \gamma_{c-1}(\mathrm{IA}(M_n))$.*

*Proof.* For $i, j \in [n]$, with $i \neq j$, we write $\pi_{ij}$ for the tame IA-automorphism of $M_n$ satisfying the conditions $\pi_{ij}(x_i) = x_i[x_i, x_j]$ and $\pi_{ij}(x_r) = x_r$, with $r \neq i$. Furthermore, we write $\sigma_{ij}$ for the tame automorphism of $M_n$ satisfying the conditions $\sigma_{ij}(x_i) = x_j$, $\sigma_{ij}(x_j) = x_i$ and $\sigma_{ij}(x_k) = x_k$, with $k \neq i, j$.

1. We shall induct on $c$. For $c = 2$, our claim is trivially true. Thus, for some $c \geq 2$, we assume that $\tau_{1,(i_1,\ldots,i_c)} \in \gamma_{c-1}(\mathrm{IA}(M_n))$ for all $i_1, \ldots, i_c \in [n] \setminus \{1\}$, with $i_1 \neq i_2$. Let $\tau_{1,(k_1,\ldots,k_{c+1})}$, with $k_1, \ldots, k_{c+1} \in [n] \setminus \{1\}$, $k_1 \neq k_2$, and consider the following two cases.

   (a) Let $k_1 \neq k_3$. It is easily checked that $\tau_{1,(k_1,k_2,\ldots,k_{c+1})} = [\pi^{-1}_{k_1,k_2}, \tau^{-1}_{1,(k_1,k_3,\ldots,k_{c+1})}]$. By our inductive hypothesis, $\tau_{1,(k_1,k_3,\ldots,k_{c+1})} \in \gamma_{c-1}(\mathrm{IA}(M_n))$. Since $\pi_{k_1,k_2} \in \mathrm{IA}(M_n)$, we have $\tau_{1,(k_1,k_2,\ldots,k_{c+1})} \in \gamma_c(\mathrm{IA}(M_n))$ for all $k_1, \ldots, k_{c+1} \in [n] \setminus \{1\}$, with $k_1 \neq k_2$.

   (b) Let $k_1 = k_3$. By direct calculations, $\tau_{1,(k_1,k_2,\ldots,k_{c+1})} = [\tau^{-1}_{1,(k_2,k_1,k_4,\ldots,k_{c+1})}, \pi^{-1}_{k_2,k_1}]$. As before, we get $\tau_{1,(k_1,k_2,\ldots,k_{c+1})} \in \gamma_c(\mathrm{IA}(M_n))$ for all $k_1, \ldots, k_{c+1} \in [n] \setminus \{1\}$, with $k_1 \neq k_2$.

2. For $c = 2$, our claim is trivially true. Thus, we may assume that $c \geq 3$. By Lemma 5 (1), we have $\tau_{123,(0,r_2,\ldots,r_n)} \in \gamma_{c-1}(\mathrm{IA}(M_n))$ for all non-negative integers $r_2, \ldots, r_n$, with $r_2 + \cdots + r_n = c-3$. We point out that $\tau_{123,(0,\ldots,0)} = \tau_{1,(2,3)} \in \mathrm{IA}(M_n)$. For a non-negative integer $s_1$, let $\psi_{1,s_1} = [\tau^{-1}_{1,(2,3)}, {}_{s_1}(\pi_{2n}\pi_{3n})^{-1}]$. Clearly, $\psi_{1,s_1} \in \gamma_{s_1+1}(\mathrm{IA}(M_n))$. By using induction on $s_1$, we show that $\psi_{1,s_1}$ satisfies the conditions $\psi_{1,s_1}(x_1) = x_1[x_2, x_3, {}_{s_1}x_n]^{-1}$ and $\psi_{1,s_1}(x_j) = x_j$, with $j \neq 1$. Since $M_n$ is metabelian, we have $[x[x,z], y[y,z]] = [x,y][x,y,z]$ for all $x, y, z \in M_n$. By direct calculations, for all $s_1 \geq 0$,
$$\tau_{123,(s_1+1,0,\ldots,0)} = [\pi_{1n}, (\psi^{-1}_{1,s_1})^{\sigma_{1n}}]. \tag{3.8}$$
Since $\psi_{1,s_1} \in \gamma_{s+1}(\mathrm{IA}(M_n))$, by the equation (3.8), $\tau_{123,(s_1+1,0,\ldots,0)} \in \gamma_{s_1+2}(\mathrm{IA}(M_n))$ for all $s_1 \geq 0$. Hence, for all $r_1 \geq 0$, we have $\tau_{123,(r_1,0,\ldots,0)} \in \gamma_{r_1+1}(\mathrm{IA}(M_n))$. For non-negative integers $s_1, \ldots, s_n$, by inductive arguments and since $M_n$ is metabelian, it may be easily verified that, for all $4 \leq j \leq n-1$,
$$\tau^{-1}_{123,(s_1+1,0,\ldots,0,s_n)} = [\tau^{-1}_{123,(s_1+1,0,\ldots,0)}, {}_{s_n}(\pi_{2n}\pi_{3n})^{-1}], \tag{3.9}$$
$$\tau^{-1}_{123,(s_1+1,0,\ldots,0,s_j,s_{j+1},\ldots,s_n)} = [\tau^{-1}_{123,(s_1+1,0,\ldots,0,s_{j+1},\ldots,s_n)}, {}_{s_j}(\pi_{2j}\pi_{3j})^{-1}], \tag{3.10}$$
$$\tau^{-1}_{123,(s_1+1,0,s_3,s_4,\ldots,s_n)} = [\tau^{-1}_{123,(s_1+1,0,0,s_4,\ldots,s_n)}, {}_{s_3}\pi^{-1}_{23}], \tag{3.11}$$

and

$$\tau^{-1}_{123,(s_1+1,s_2,s_3,\ldots,s_n)} = [\tau^{-1}_{123,(s_1+1,0,s_3,\ldots,s_n)}, {}_{s_2}\pi^{-1}_{32}]. \tag{3.12}$$

Thus, for non-negative integers $r_1, r_2, \ldots, r_n$, with $r_1 \geq 1$ and $r_1+r_2+\cdots+r_n = c-2$, by the equations (3.9)-(3.12), we get

$$\tau^{-1}_{123,(r_1,r_2,\ldots,r_n)} = [\tau^{-1}_{123,(r_1,0,\ldots,0)}, {}_{r_n}(\pi_{2n}\pi_{3n})^{-1}, \ldots, {}_{r_4}(\pi_{24}\pi_{34})^{-1}, {}_{r_3}\pi^{-1}_{23}, {}_{r_2}\pi^{-1}_{32}]. \tag{3.13}$$

Since, as shown above, $\tau_{123,(r_1,0,\ldots,0)} \in \gamma_{r_1+1}(\mathrm{IA}(M_n))$ for all $r_1 \geq 0$, by the equation (3.13), we have the required result.

# 5 Proof of Theorem 1

For $n = 2$, Theorem 1 (1) follows from Corollary 1 and, for $n = 3$, Theorem 1 (2) follows from Proposition 1. For $n \geq 4$ and $c \geq 4$, by Proposition 2, the $\mathbb{Z}\mathrm{GL}_n(\mathbb{Z})$-module $\mathcal{L}^c(\mathrm{IA}(M_n))$ is generated by $\overline{\tau}_{123,(r_1,\ldots,r_n)}$ and $\overline{\xi}_{v(k_1,\ldots,k_{c-1})}$. Since $\xi_{v(k_1,\ldots,k_{c-1})}$ is an inner automorphism of $M_n$ and since $Z(M_n/\gamma_{c+1}(M_n)) = \gamma_c(M_n)/\gamma_{c+1}(M_n)$ for all $c$ (see, for example, [14, Chapter 3, Section 6]), it follows from Lemma 3 that $\xi_{v(k_1,\ldots,k_{c-1})} \in \gamma_{c-1}(\mathrm{IA}(M_n))$. Furthermore, by Lemma 5 (2), we have $\tau_{123,(r_1,r_2,\ldots,r_n)} \in \gamma_{c-1}(\mathrm{IA}(M_n))$ and so, we obtain the desired result.